\documentclass[10pt,a4paper]{article}
\usepackage[latin1]{inputenc}
\usepackage{amsmath}
\usepackage{amsfonts}
\usepackage{amssymb}
\usepackage{makeidx}
\usepackage{MnSymbol}
\usepackage{hyperref}
\usepackage{amsthm}
\author{Sol Schwartzman\\
c/o Math Department\\
University of Rhode Island\\
Kingston, RI\\
solschwartzman@gmail.com}
\title{Zeroes of Real Valued Eigenfunctions}

\newtheorem*{theorem*}{Theorem}
\newtheorem*{lemma*}{Lemma}
\newtheorem*{corollary*}{Corollary}

\begin{document}
\maketitle

\normalsize
Suppose we are given a symmetric operator $T$ acting on a subspace of $L_2(M_n,\mu)$ where $M_n$ is an $n$-dimensional metrizable connected manifold and $\mu$ is a measure that is positive on open sets in $M_n$.  Then there is at most one eigenspace that contains a real valued eigenfunction whose nodal set has dimension less than $n-1$.

\begin{center}Spectral Theory (math.ST)\end{center}

\clearpage
Let $M_n$ be a connected $n$-dimensional metrizable manifold and let $\mu$ be a measure on $M_n$ that is positive on open sets.  Suppose that $T$ is a symmetric operator acting on a subspace of $L_2(M_n,\mu)$.

\begin{theorem*}There is at most one eigenspace that contains a real valued continuous function $\psi$ for which the nodal set $N_\psi$ of its zeroes has dimension less than $n-1$.
\end{theorem*}

\textbf{Proof:}  First we will need the following.

\begin{lemma*}  If $\psi$ is any continuous real valued function on $M_n$ whose nodal set $N_\psi$ has dimension less than $n-1$, then $\psi$ is of constant sign on the complement of $N_\psi$
\end{lemma*}
\textbf{Proof of lemma:}  The complement of $N_\psi$ is connected [1] and $\psi(x) \neq 0$ for any $x$ outside $N_\psi$.  This establishes the lemma.

Next, let $\psi_1$ and $\psi_2$ be real valued continuous eigenfunctions corresponding to different eigenvalues, and suppose that the dimensions of $N_{\psi_1}$ and $N_{\psi_2}$ are each less than $n-1$.  Since $N_{\psi_1 \psi_2}$=$N_{\psi_1} \cup N_{\psi_2}$, it follows from the sum theorem for dimension that the dimension of $N_{\psi_1 \psi_2}$ is less than $n-1$.  Our lemma implies that either $\psi_1 \psi_2$ or $-\psi_1 \psi_2$ is positive outside  $N_{\psi_1 \psi_2}$.  It follows that the integral of $\psi_1 \psi_2$ over $M_n$ cannot be zero, which gives us a contradiction and establishes our theorem.

$\blacksquare$ 

We will say that an eigenvalue $\lambda$ for which there is a $\psi$ in the corresponding eigenspace whose nodal set $N_\psi$ is of dimension less than $n-1$ is an \emph{exceptional eigenvalue}.

\begin{corollary*}
If $M_n$ and $\psi$ are real analytic and the eigenvalue for $\psi$ is not exceptional, then the dimension of $N_\psi$ must equal $n-1$.
\end{corollary*}

\textbf{Proof:}  If the dimension of $N_\psi$ equalled $n$, then $N_\psi$ would have to contain an open set[1].

$\blacksquare$

\end{document}